\documentclass[a4paper,11pt]{amsart}

\usepackage[english]{babel}
\usepackage{amsmath, enumerate, amsfonts, amssymb, amsthm}
\usepackage[dvips]{graphicx}

\newtheorem{defin}{Def\mbox{}inition}[section]
\newtheorem{theo}[defin]{Theorem}
\newtheorem{prop}[defin]{Proposition}
\newtheorem{lem}[defin]{Lemma}
\newtheorem{cor}[defin]{Corollary}
\newtheorem{rem}[defin]{Remark}

\newtheorem*{defin*}{Def\mbox{}inition}
\newtheorem*{theo*}{Theorem}
\newtheorem*{prop*}{Proposition}
\newtheorem*{lem*}{Lemma}
\newtheorem*{cor*}{Corollary}
\newtheorem*{rem*}{Remark}
\newtheorem*{aff*}{Claim}
\newtheorem*{affs*}{Claims}
\newtheorem*{nota*}{Notation}

\newtheorem*{theobis*}{Theorem 2 bis}

\newtheorem{theoi}{Theorem}

\newcommand{\dps}{\displaystyle}

\newcommand{\R}{\mathbb{R}}
\newcommand{\C}{\mathbb{C}}
\newcommand{\N}{\mathbb{N}}
\newcommand{\Z}{\mathbb{Z}}
\newcommand{\Q}{\mathbb{Q}}

\newcommand{\An} {\mathbf{A}_n}

\newcommand{\D}{\mathcal{D}}
\newcommand{\Dn}{\mathcal{D}_n}
\newcommand{\Dnp}{\mathcal{D}_{n+p}}

\newcommand{\Dnpz} {\mathcal{D}_{n+p}\langle z \rangle}
\newcommand{\Dm}{\mathcal{D}_m}
\newcommand{\Dmz}{\mathcal{D}_m\langle z \rangle}
\newcommand{\Dnd}{\mathcal{D}_{n+2}}

\newcommand{\dxi}{\partial _{x_i}}

\newcommand{\dx}[1]{\partial _{x_{#1}}}
\newcommand{\dt}[1]{\partial _{t_{#1}}}
\newcommand{\ddx}{\partial _x}
\newcommand{\ddt}{\partial _t}

\newcommand{\dxsur}[2]{\frac{\partial {#1}}{\partial x_{#2}}}

\newcommand{\B}{\mathcal{B}}

\newcommand{\ord}{\mathrm{ord}}
\newcommand{\gr}{\mathrm{gr}}
\renewcommand{\mp}{\mathrm{lm}} 

\newcommand{\DN}{\mathcal{N}} 
\newcommand{\Exp}{\mathrm{Exp}}

\newcommand{\Vb}{\overline{V}}
\newcommand{\U}{\mathcal{U}}
\newcommand{\UV}{\mathcal{U}_V}
\newcommand{\E}{\mathcal{E}}
\newcommand{\EV}{\mathcal{E}_V}
\renewcommand{\L}{\mathcal{L}}

\begin{document}

\title[Existence of Bernstein-Sato polynomials]{Existence of Bernstein-Sato polynomials by using the analytic Gr\"obner fan}

\author{Rouchdi BAHLOUL}
\email{rouchdi@math.kobe-u.ac.jp}
\address{Department of Mathematics\\Faculty of Science\\Kobe University\\1-1, Rokkodai,\\Nada-ku\\Kobe 657-8501\\Japan}




\begin{abstract}
In 1987, C.~Sabbah proved the existence of Bernstein-Sato polynomials
associated with several analytic functions. The purpose of this
article is to give a more elementary and constructive proof of the
result of C.~Sabbah based on the notion of the analytic Gr\"obner fan
of a $\D$-module.
\end{abstract}

\maketitle




This paper is a translation of \cite{comp}.

\section*{Introduction and statement of the main results}\label{sec:intro}

Fix two integers $n \ge 1$ and $p\ge 1$ and $v \in \N^p \smallsetminus \{0\}$. Let $x=(x_1,\ldots,x_n)$ and $s=(s_1,\ldots,s_p)$ be two systems of variables. Consider $f_1,\ldots, f_p \in \C\{x\}=\C\{x_1,\ldots, x_n\}$ and denote by $\Dn$ the ring of differential operators with coefficients in $\C\{x\}$.
For $b(s) \in \C[s]=\C[s_1,\ldots,s_p]$, consider the following identity: 
\[(\star) \qquad b(s) f^s \in \D_n[s] f^{s+v},\]
where $f^{s+v}=f_1^{s_1+v_1} \cdots f_p^{s_p+v_p}$. A polynomial
$b(s)$ satisfying such an identity is called a Bernstein-Sato polynomial (associated with $f=(f_1,\ldots,f_p)$). The set of these polynomials form an ideal called the Bernstein-Sato ideal and denoted by $\B^v(f)$. Let us mention that usually, $v$ is taken as $(1,\ldots, 1)$ or $(0, \ldots, 0, 1, 0, \ldots,0)$ where $1$ is in the $j$-th position, $j \in \{1,\ldots, p\}$.\\
Let us give some historical recalls. In the case where $p=1$ and $f$ is a polynomial, I.N. Bernstein \cite{bernstein} showed that the ideal $\B^v(f)$ is not zero (in this case, in $(\star)$, $\Dn$ is replaced with the Weyl algebra $\An(\C)$, i.e. the ring of differential operators with polynomial coefficients). Again for $p=1$ but in the analytic case, the fact that $\B^v(f)$ is not zero was proved by J.~E.~Bj\"ork \cite{bjork} with similar methods as that of \cite{bernstein}. In the same case, let us cite M.~Kashiwara \cite{kashiwara} who published another proof and showed moreover that the unitary generator of the Bernstein-Sato ideal has rational roots.
Now, for $p\ge 2$ the proof in the polynomial case is an easy generalization of that by I.~N.~Bernstein, which can be found in \cite{lichtin}. In the analytic case, the proof of the existence of a non zero Bernstein-Sato polynomial was given by C.~Sabbah (\cite{sabbah1} and \cite{sabbah2}). Let us cite the contribution of A.~Gyoja \cite{gyoja} who showed moreover that $\B^v(f)$ contains a rational non zero element.

The goal of the present paper is a development of the proof by C.~Sabbah. More precisely, we can decompose the proof of C.~Sabbah into two main steps: the first one uses arguments similar to that used by M.~Kashiwara in the case $p=1$, the second one essentially consists in a finiteness result that reduces the problem to the first step. The second step of the proof by C.~Sabbah is based on an ``adapted fan''. The existence of such a fan is done in (\cite{sabbah1} theor. A.1.1). Unfortunately, there is a gap in the proof of theorem A.1.1. (see also the comments after th.~S1 in the present paper). In this paper, we propose a more elementary and constructive statement and proof of the second step, by avoiding the sensitive notion of ``adapted fan''.

In order to motivate the reading of this article, we recall (without all the details) the proof by C.~Sabbah. In this recall, we will emphasize the result (th. S1) for which we will give a more constructive statement. Let us mention that most of the notions or notations introduced below and useful for the sequel will be detailed in the next sections.

We denote by $\Dnp$ the ring of differential operators with coefficients in $\C\{x,t\}=\C\{x_1,\ldots,x_n, t_1,\ldots,t_p\}$. Following the method of B.~Malgrange \cite{malgrange}, we make $\Dnp$ act on $\C\{x\}[\frac{1}{f_1\cdots f_p}, s] f^s$. We denote by $I$ the (left\footnote{Throughout the text, ideal shall mean left ideal}) ideal annihilator of $f^s$ in $\Dnp$ and $M$ the quotient $M=\Dnp/I$.

For $j=1,\ldots,p$, denote by $V_j(\Dnp)$ the Kashiwara-Malgrange $V$-filtration associated with the variable $t_j$ et set $V=(V_1,\ldots, V_p)$ the (multi)filtration of $\Dnp$ indexed by $\Z^p$: for $w\in \Z^p$,
\[V_w(\Dnp)=\bigcap_{j=1}^p \{V_j\}_{w_j}(\Dnp).\]
It induces a filtration $V(M)$ on $M$ where for $w \in \Z^p$, $V_w(M)$ is the image of $V_w(\Dnp)$ by the projection $\Dnp \to M=\Dnp/I$.

For $j=1,\ldots,p$, we identify the filtration $V_j$ with the linear form on $\N^{2n+2p}$ given by $V_j(\alpha,\mu, \beta,\nu)=\nu_j -\mu_j$ (where $\alpha,\beta \in \N^n$, $\mu,\nu \in \N^p$ and $\alpha, \mu, \beta, \nu$ corresponds respectively to $x, t, \ddx, \ddt$).
Denote by $\UV=\sum_{j=1}^p \R_{\ge 0} V_j$. We identify $\UV$ with $(\R_{\ge 0})^p$. Each $L$ in $\UV \cap \N^p$ (i.e. $L$ with integral coefficients) gives rise to a natural filtration $V^L$ on $\Dnp$ and $M$ indexed by $\Z^p$ given by:
\[V^L_k(M) = \sum_{\{w\in \Z^p;\, L(w)\le k\} } V_w(M),\]
where $L(w)=l_1 w_1+\cdots+l_p w_p$ if $L=(l_1,\ldots,l_p) \in \N^p$.\\
\\
Now, let us recall the two main steps of the proof by C.~Sabbah.
\begin{description}
\item[Step 1]
\begin{theo*}(\cite{sabbah1} th. 3.1.1, see also \cite{gyoja} 2.9 and
  2.10)
For any $L \in \UV \cap \N^p$, there exists a non zero polynomial $b\in \C[\lambda]$ in one variable such that for any $k\in \Z$,
\[b\big(L(-\dt{1} t_1, \ldots, -\dt{p}t_p)-k\big) V^L_k(M) \subset
V^L_{k-1}(M).\]
\end{theo*}
The proof of this theorem uses arguments similar to that used by M.~Kashiwara \cite{kashiwara} in the case $p=1$.

Denote by $b_L$ the monic polynomial with minimal degree satisfying the previous identity. The contribution of A.~Gyoja \cite{gyoja} consists in the fact that $b_L$ has roots in $\Q_{<0}$.

\item[Step 2]

We are going to introduce two other filtratrions on $M$.
\begin{enumerate}
\item
Let $\sigma$ be a convex rational cone in $\R_{\ge 0}^p$. We denote by $\L(\sigma)$ the set of the primitive elements of the $1$-skeleton of $\sigma$ (i.e. $L\in \L(\sigma)$ if and only if the line generated by $L$ is in the $1$-skeleton of $\sigma$ and the coefficients of $L$ are integral and without a common factor $\ne 1$). For any $w\in Z^p$, we set:
\[{}^\sigma V_w(M)= \sum_{\{w'\in \Z^p \, | \forall L \in \L(\sigma)
  \, L(w')\le L(w)\}} V_{w'}(M).\]
(See figure \ref{fig:Vsigma} where $p=2$ and $\sigma$ is generated by $L_1$ and $L_2$: $m\in {}^\sigma V_w(M) \iff m$ is represented by an operator $P\in \Dnp$ which has a Newton diagram included in the cross-ruling.)
\begin{figure}[h!]
\centering
\includegraphics[angle=0, width=11cm]{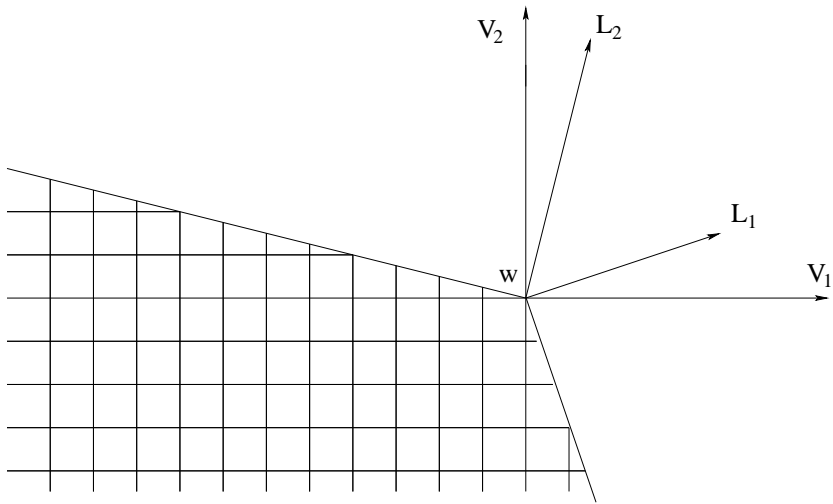}
\caption{${}^\sigma V_w(M)$}\label{fig:Vsigma}
\end{figure}

For any $w\in \Z^p$, it is easy to see that we have the following inclusions:
\[{}^\sigma V_w(M) \subseteq \bigcap_{L\in \sigma} V^L_{L(w)}(M)
\subseteq \bigcap_{L\in \L(\sigma)} V^L_{L(w)}(M),\]
the first one being trivial and the second one resulting from the fact that $\L(\sigma)$ is included in the closure of $\sigma$.

\begin{theo*}[S1](\cite{sabbah1} th. A.1.1 and prop. 2.2.1)
There exists a fan $\Sigma$ (called fan adapted to $V(M)$) made of convex rational polyhedral cones such that for any cone $\sigma\in \Sigma$ and any $w\in \Z^p$, we have:
\[{}^\sigma V_w(M) = \bigcap_{L\in \L(\sigma)} V^L_{L(w)}(M).\]
\end{theo*}
This is the theorem for which we will give a more elementary statement and proof. Let us say a word about the proof. In the appendix of \cite{sabbah1}, in collaboration with F.J.~Castro Jim\'enez, C.~Sabbah states the existence of a fan $\Sigma$ adapted to the filtration $V(M)$ (th.~A.1.1) then he shows, in prop. 2.2.1, that for any cone in such a fan, we have the previous equality. The exists a gap in the proof of theorem A.1.1. Indeed, the proof is based on a division with parameters which gives rise to formal power series in the derivation variables $\dxi$.
\item
Let $\Vb(M)$ be the filtration indexed by $\Z^p$ defined by:
\[\Vb_w(M)=\bigcap_{L\in \UV} V^L_{L(w)}(M).\]
A consequence of theorem~S1 is

\begin{cor*}[S2]\
\begin{enumerate}
\item
For any $w\in \Z^p$,
\[\Vb_w(M)= \bigcap_{L\in \L(\Sigma)} V^L_{L(w)}(M),\]
where $\L(\sigma)$ denotes the set of the primitive elements of the $1$-skeleton of $\Sigma$.
\item
There exists $\kappa \in \N^p$ such that for any $w\in \Z^p$,
\[V_w(M) \subset \Vb_w(M) \subset V_{w+\kappa}(M).\]
\end{enumerate}
\end{cor*}

Let us make some comments about the proof of this result. Assertion (a) of S2 trivially follows from S1, indeed:
\begin{eqnarray*}
\Vb_w(M) & = & \bigcap_{\sigma \in \Sigma} \big(\bigcap_{L\in \sigma}
V^L_{L(w)}(M) \big)\\
 & = & \bigcap_{\sigma \in \Sigma} \big( \bigcap_{L\in \L(\sigma)}
V^L_{L(w)}(M)\big) \text{ by S1}\\
 & = & \bigcap_{L\in \L(\Sigma)} V^L_{L(w)}(M).
\end{eqnarray*}
Concerning assertion (b), C.~Sabbah proves that if a fan $\Sigma$ satisfies the satement of theorem S1 then $\Vb(M)$ is a good $V(\Dnp)$ filtration, thus by a usual result (on the comparison of good filtrations), we obtain the wanted inclusions (remark that the first inclusion is trivial and it is the second one that we are interested in).
\end{enumerate}
Let us end the recall of the proof. Denote by $\delta$ the class of $1$ in the quotient $M=\Dnp/I$. Put
\[b(s)= \prod_{L\in \L(\Sigma)} \Big(\prod_{-L(v+\kappa)<k\le 0}
b_L\big(L(s)-k\big)\Big).\]
By assertions (a) and (b) of S2, we obtain:
\[b(-\ddt t) \delta \in \Vb_{-v-\kappa}(M) \subset V_{-v}(M),\]
which means that $b(s) \in \B^v(f)$.

\end{description}

Let us now state the main results of the present paper.

Let $I$ be an ideal in $\Dm$ (ring of differential operators with analytic coefficients in $m$ variables). In \cite{acg}, A.~Assi, F.J.~Castro Jim\'enez and M.~Granger introduced  the set $\U$ of the linear forms $L$ for which the naturally associated filtration on $\Dm$ is compatible with the non commutative structure of $\Dm$ (see paragraph 1.2). They studied the behaviour of the graded ring $\gr^L(I)$ when $L$ moves in $\U$. Consider the relation on $\U$ such that $L$ is in relation with $L'$ if the graded rings $\gr^L(h(I))$ and $\gr^{L'}(h(I))$ are equal (here, $h(I)$ denotes the homogenization of $I$, we shall recall it later). This relation, which is an equivalence relation, gives rise to a partition of $\U$ made of convex rational polyhedral cones. This partition is called the (analytic) Gr\"obner fan associated with $h(I)$ et denoted by $\E(h(I))$.

Now, put $m=n+p$ and resume the previous notations. As we will see, we can naturally include $\UV$ into $\U$. Denote by $\EV=\EV(h(I))$ the fan on $\UV$ obtained by restriction of $\E(h(I))$. Here are the results we aim to prove.

\begin{theoi}\label{theo}
For any cone $\sigma$ in $\EV$,
\[{}^\sigma V_w(M) = \bigcap_{L\in \L(\sigma)} V^L_{L(w)}(M).\]
\end{theoi}
The counterpart of corollary S2 follows by replacing the $1$-skeleton of $\Sigma$ with that of $\EV$. Indeed, the proof of corollary S2 works for any fan from the moment that it satisfies the statement of theorem S1.

\begin{theoi}\label{theo2}
For $p=2$, there exists $\kappa \in \N^2$ that we can compute from the Gr\"obner bases associated with each cone of $\EV$ (this computation shall be detailed in section 3) such that for any $w\in \Z^2$,
\[\Vb_{w+\kappa}(M) \subset V_w(M).\]
\end{theoi}
For $p\ge 3$, a generalization of this result seems to give rise to technical difficulties hard to solve.

Finally, here are the main contributions of the present paper.
\begin{itemize}
\item
A constructive statement and proof of the key theorem S1 which provides a more elementary and constructive approach of step 2 in the proof by C.~Sabbah and which avoids the notion of adapted fan.
\item
For $p=2$, a completely constructive proof of step 2 in the proof of C.~Sabbah.
\end{itemize}
I point out that this article is a part of my thesis \cite{bahloul} in which one can found another proof of the counterpart (i.e. with $\L(\EV)$ instead of $\L(\Sigma)$) of statement (a) of S2, without using theorem \ref{theo}.\\
\\
In a first section, we shall make some recalls concerning the division theorem in $\Dm$ as in \cite{acg}, the notions of standard basis and Gr\"obner fan. In section 2, we shall give the proof of theorem \ref{theo} and in section 3 the proof of theorem \ref{theo2}.

\section{Recalls and preparatory results}\label{sec:rappels}

In the following paragraphs, we shall recall without giving proofs some notions and results which are helpful for the sequel.

\subsection{Homogenization}

As \cite{castro-narvaez} in the algebraic case, let us introduce the homogenized ring $\Dmz$ with which the authors of \cite{acg} introduced the analytic Gr\"obner fan.

In this paragraph and until 1.3.2, $x=(x_1,\ldots,x_m)$ and $\Dm$ denotes the ring of differential operators with coefficients in $\C\{x\}$. We define $\Dmz$ as the $\C\{x\}$-algebra generated by $\dx{1}, \ldots, \dx{m}, z$ where the only non trivial commutation relation are:
\[[\dx{i}, c(x)]=\dxsur{c(x)}{i} z \text{ for } i=1,\ldots,n \text{
 and } c(x) \in \C\{x\}.\]
In $\Dm$, we denote by $\deg(P)$ the total degree of $P$ in the $\dxi$. Consider the associated filtration. We can extend it to $\Dmz$ by considering the total degree in the $\dxi$ and $z$. This filtration provides to $\Dmz$ a structure of graded algebra:
\[\Dmz=\bigoplus_{d\in \N} \Dmz_d\quad  \text{ with }
\quad \Dmz_d=\!\!\!\!\bigoplus_{k+|\beta|=d}
\C\{x\}\ddx^\beta z^k,\]
where $\beta\in \N^m$, $\ddx^\beta=\dx{1}^{\beta_1} \cdots
\dx{m}^{\beta_m}$ and $|\beta|=\beta_1+\cdots+\beta_m$. Remark that the filtration of $\Dm$ given by $\deg$ gives rise to a Rees algebra which is isomorphic to $\Dmz$.\\
We say that an operator $P\in \Dmz$ is homogeneous (of degree $d$) if $P \in \Dmz_d$.\\
For $P\in \Dm$, we define its homogenization $h(P)\in \Dmz$ as follows.\\
Write $P=\sum_{\beta}c_\beta(x) \ddx^\beta$ and set $h(P)=\sum_\beta c_\beta(x)\ddx^\beta z^{d-|\beta|}$ where $d=\deg(P)$, thus $h(P)$ is homogeneous of degree $\deg(P)$.\\
Now, for an ideal $I$ in $\Dm$, we define $h(I)$ as the ideal of $\Dmz$ generated by the set of $h(P)$ for $P\in I$.

\subsection{Filtrations, divisions and standard bases}

We are going to recall the notions of a filtration on $\Dm$ and $\Dmz$ and the division theorem in $\Dmz$ as in \cite{acg}. We shall also recall the notion of standard basis and minimal reduced standard basis.

Let $\U$ be the set of the linear forms $L:\R^{2m} \to \R$,
$L(\alpha,\beta)=\sum_1^n e_i \alpha_i +\sum_1^n f_i \beta_i$ where for any $i=1,\ldots,n$, $e_i\le 0$ and $e_i+f_i\ge 0$. We extend $\U$ to $\R^{2m+1}$ by seting $L(\alpha,\beta,k)=L(\alpha,\beta)$.\\
For $P\in \Dmz$ (resp. $P\in \Dm$), write $P=
\sum_{\alpha,\beta,k} a_{\alpha, \beta, k} x^\alpha \ddx^\beta z^k$
(resp. with $a_{\alpha, \beta, k}=0$ for $k >0$ if $P\in \Dm$). We define the Newton diagram of $P$, denoted by $\DN(P)$ in $\N^{2m+1}$ (resp. in $\N^{2m}$), as the set of $(\alpha,\beta,k) \in \N^{2m+1}$ for which $a_{\alpha,\beta,k} \ne
0$.\\
Given $L\in \U$ and $P\in \Dmz$ (or $P\in \Dm$), we define $L$-order of $P$ as $\ord^L(P)=\max L(\DN(P))$. This order induces a filtration $F^L$ on $\Dmz$ or $\Dm$ indexed by $L(\N^{2m+1})$ given by:
\[F^L_k(\Dmz)=\{P\in \Dmz; \, \ord^L(P) \le k\}\]
and an associated graded ring $\dps\gr^L(\Dmz)= \bigoplus_{k\in
L(\N^{2m+1})} F^L_k(\Dmz)/F^L_{<k}(\Dmz)$.\\
For $P\in \Dmz$, we denote by $\sigma^L(P)$ the principal symbol of 
$P$, i.e. the class of $P$ in the quotient $F^L_k(\Dmz)/F^L_{<k}(\Dmz)$ where
$k=\ord^L(P)$. If $J$ is an ideal in $\Dmz$, a filtration $F^L(J)$ is induced and gives rise to a graded ideal $\gr^L(J)$ of $\gr^L(\Dmz)$, which is generated by the set of the $\sigma^L(P)$ for $P\in J$.\\
\\
For a form $L\in \U$, we define two orders: $<_L$ on $\N^{2m}$ and $<_L^h$ on $\N^{2m+1}$:
\[(\alpha,\beta) <_L (\alpha',\beta') \iff 
\begin{cases}
L(\alpha,\beta) < L(\alpha',\beta') \\
\text{or } \big(  \text{ equality and } |\beta|<|\beta'| \big) \\
\text{or } \big( \text{ equality and } (\alpha,\beta) >_0
(\alpha',\beta') \big),
\end{cases}\]
where $<_0$ is a well order compatible with sums, which is fixed for the sequel;

\[(\alpha,\beta,k) <_L^h (\alpha',\beta',k') \iff
\begin{cases}
k+|\beta|<k'+|\beta'| \\
\text{or } \big( \text{ equality and } (\alpha,\beta) <_L
(\alpha',\beta') \big).
\end{cases}\]

For $P\in \Dmz$, we define the leading exponent $\exp_{<_L^h}(P)= \max _{<_L^h}(\DN(P))$ and the leading monomial $\mp_{<_L^h}(P)=(x,\ddx,z)^{\exp_{<_L^h}(P)}$. We do the same for $P\in \Dm$ and with $<_L$ instead of $<_L^h$. Note that $\exp_\prec(PQ)=\exp_\prec(P) +\exp_\prec(Q)$ if $\prec$ is compatible with sums, which is the case for $<_L$ and $<_L^h$. Let us recall the division theorem in $\Dmz$ given in \cite{acg}.

Let $L\in\U$. Let $Q_1, \ldots, Q_r$ be a family of operators in $\Dm$. Denote by $(\Delta_1,\ldots, \Delta_r,\bar{\Delta})$ the partition of $\N^{2m+1}$ defined from the $\exp_{<_L^h}(Q_j)$:
\begin{itemize}
\item $\Delta_1=\exp_{<_L^h}(Q_1)+\N^{2m+1}$
\item $\Delta_j=(\exp_{<_L^h}(Q_j)+\N^{2m+1}) \smallsetminus
  (\bigcup_{i=1}^{i=j-1} \Delta_i)$ for $j=2,\ldots,r$
\item $\bar{\Delta}=\N^{2m+1} \smallsetminus (\bigcup_{j=1}^{j=r}
  \Delta_j)$
\end{itemize}
\begin{theo} (\cite{acg} Th. 7)\label{theo:div}
For any $P\in \Dmz$, there exists a unique $(q_1,\ldots,q_r,\\ R) \in
(\Dmz)^{r+1}$ such that:
\begin{enumerate}
\item $P=q_1 Q_1 + \cdots + q_r Q_r+R$
\item for any $j=1,\ldots,r$, if $q_j \ne 0$ then $\DN(q_j)+\exp_{<_L^h}(Q_j) \subset \Delta_j$
\item if $R \ne 0$ then $\DN(R) \subset \bar{\Delta}$.
\end{enumerate}
We call $R$ the remainder of the division of $P$ by the $Q_j$ w.r.t. $<_L^h$.
\end{theo}

\begin{cor}\
\begin{itemize}
\item
$\exp_{<_L^h}(P)=\max_{<_L^h}\{\exp_{<_L^h}(q_jQ_j),\, j=1,\ldots,r;
\, \exp_{<_L^h}(R)\}$.
\item
$\ord^L(P)=\max\{\ord^L(q_jQ_j),\, j=1,\ldots,r;\, \ord^L(R)\}$.
\end{itemize}
\end{cor}
Let $J$ be an ideal in $\Dmz$ and $Q_1,\ldots,Q_r \in J$. We say that $Q_1,\ldots,Q_r$ form a $<_L^h$-standard bases of $J$ if for any $P$ in $J$, the remainder of the division of $P$ by the $Q_j$ is zero. Consider the set of exponents of $J$: $\Exp_{<_L^h}(J)=\{\exp_{<_L^h}(P), \, P\in J\smallsetminus 0\}$. Given $Q_1,\ldots, Q_r \in J$, these two claims are equivalent (by the division theorem):
\begin{itemize}
\item The $Q_j$ form a $<_L^h$-standard basis of $J$.
\item $\Exp_{<_L^h}(J)=\bigcup_{j=1}^r (\exp_{<_L^h}(Q_j)+\N^{2m+1})$.
\end{itemize}
The existence of a standard basis holds by Dickson lemma which asserts that if a subset $E$ of $\N^q$ satisfies $E=E+\N^q$ (which is the case for $\Exp_{<_L^h}(J)$) then there exists $F\subset E$ finite such that $E=\cup_{e\in F} (e+\N^q)$.
It is easy to see that a standard basis is not unique in general, that is why there exists the notion of minimal reduced standard basis:
\begin{defin*}
Let $Q_1,\ldots,Q_r$ be a $<_L^h$-standard basis of $J\subset \Dmz$ and let $e_j=\exp_{<_L^h}(Q_j)$ for $j=1,\ldots,r$.
\begin{itemize}
\item
We says that it is minimal if for any finite subset $F$ of $\N^{2m+1}$, the following implication holds
\[\Exp_{<_L^h}(J)=\bigcup_{e\in F}(e+\N^{2m+1}) \Rightarrow
\{e_1,\ldots,e_r\} \subseteq F.\]
\item
We say that it is reduced if the $Q_j$ are unitary (i.e. the coefficient corresponding to the leading monomial is $1$) and if for any $j$,
\[(\DN(Q_j) \smallsetminus e_j) \subset (\N^{2m+1} \smallsetminus
\Exp_{<_L^h}(J)).\]
\end{itemize}
There exists a unique $<_L^h$-minimal reduced standard basis of an ideal in $\Dmz$.
\end{defin*}
Remark that if the ideal $J$ is homogemeous then the minimal reduced standard basis will be made of homogeneous elements. Let us end this paragraph by a result that we will apply in the next section.

\begin{lem}\label{lem:utile}
Let $<_1$ and $<_2$ be two orders on $\N^{2m+1}$ which allow one to make divisions in $\Dmz$ (for example $<_i=<_{L_i}^h$ with $L_1$ and $L_2$ being two forms). Let $\{Q_1,\ldots,Q_r\}$ be the minimal reduced standard basis of an ideal $J \subset \Dmz$ w.r.t. $<_1$. Suppose that for any $j$, $\exp_{<_1}(Q_j)$ and
$\exp_{<_2}(Q_j)$ are equal. Then $\{Q_1,\ldots,Q_r\}$ is also the minimal reduced standard basis of $J$ w.r.t. $<_2$.
\end{lem}
We omit the proof. Let us just say that what is important in a division is the set of the $\exp_<(Q_j)$ and not the order itself. Thus a division w.r.t. $<_1$ or to $<_2$ will give the same quotients and the same remainder.


\subsection{Gr\"obner fan}

In this paragraph, we shall recall the main result of \cite{acg} which describes the analytic Gr\"obner fan. We shall also introduce the $V$-Gr\"obner fan $\EV(h(I))$ which will be the main object of sections $2$ and $3$.

\subsubsection{}

Let $I$ be an ideal in $\Dm$. Consider $h(I)\subset \Dmz$ its homogenization. For $L$ and $L'$ in $\U$, we define the relation:
\[L \sim L' \iff \gr^L(\Dmz) = \gr^{L'}(\Dmz) \text{ and }
\gr^L(h(I)) = \gr^{L'}(h(I)).\]
This an equivalence relation on $\U$.

\begin{theo}\label{theo:event_Grob}\cite{acg}
The partition of $\U$ given by this relation is made of convex polyhedral rational cones. This partition denoted $\E(h(I))$ is called the (analytic) Gr\"obner fan of $h(I)$.\\
Moreover for any cone $\sigma \in \E(h(I))$, there exists 
$Q_1,\ldots,Q_r \in h(I)$ homogeneous such that:
\begin{itemize}
\item
for any $L,L'$ in $\sigma$, $\sigma^L(Q_j)=\sigma^{L'}(Q_j)$ and $\exp_{<_L^h}(Q_j)=\exp_{<_{L'}^h}(Q_j)$ for any $j$.
\item
for any $L\in \sigma$, the set $\{Q_1,\ldots,Q_r\}$ is the minimal reduced standard basis of $h(I)$ w.r.t. $<_L^h$.
\end{itemize}
\end{theo}
By the second claim, we can see that on a cone $\sigma$, the set $\Exp_{<_L^h}(h(I))$ of exponents of $h(I)$ w.r.t. $<_L^h$ is constant when $L$ runs over $\sigma$.\\

\subsubsection{}

From now on, we shall work in $\Dnp$ and $\Dnpz$, that is $m=n+p$. For $j=1,\ldots,p$, we denote by $V_j \in \U$ the linear form given by: $V_j(\alpha,\mu,\beta,\nu)=\nu_j -\mu_j$ where $\alpha,\beta \in \N^n$ and $\mu,\nu\in \N^p$. This form gives rise to a filtration that we also denote by $V_j$ and which is nothing but the Kashiwara-Malgrange $V$-filtration associated with the variable $t_j$ (recall that in $\Dnp$, the variables are $x=(x_1,\ldots,x_n)$ and $t=(t_1,\ldots, t_p)$). We denote by $V$ the multifiltration $V=(V_1,\ldots,V_p)$.\\
We set $\UV\subset \U$ to be the subset of the linear forms $L$ of the form:
\[L=l_1 V_1+\cdots+l_p V_p,\]
with $(l_1,\ldots,l_p)\in (\R_{\ge 0})^p$, thus we shall identify $\UV$ and $(\R_{\ge 0})^p$.

From now on, for $L\in \UV$, we shall denote by $V^L$ the filtration associated with $L$ (following the notations of \cite{sabbah1}). Remark that for $L\in \UV$, we can identify $\gr^L(\Dnpz)$ to a subring of $\Dnp$. Indeed, after reordering the variables, assume that $L=l_1 V_1+\cdots +l_e V_e$ with $0\le e \le p$ and none of the $l_j$ is zero (by convention, if $e=0$ then $L=0$ and $\gr^L(\Dnpz)= \Dnpz$). In this case, we have
\[\gr^L(\Dnpz)= \C\{x_{e+1}, \ldots, x_p\}[x_1, \ldots, x_e][\dx{1}, \ldots, \dx{n}]\langle z\rangle,\]
in which the commutation relations coincides with those of $\Dnpz$. Thus we will consider that all the calculations are made in $\Dnpz$. Now, consider the restriction of $\E(h(I))$ to $\UV$. We obtain the $V$-Gr\"obner fan of $h(I)$, denoted $\EV(h(I))$. Similarly to \ref{theo:event_Grob}, we obtain:

\begin{cor}\label{cor:V_event_Grob}
For any cone $\sigma$ of $\EV(h(I))$, there exists $Q_1,\ldots,Q_r \in h(I)$ homogeneous such that:
\begin{itemize}
\item for any $L,L' \in \sigma$, $\sigma^L(Q_j)=\sigma^{L'}(Q_j)$
  and $\exp_{<_L^h}(Q_j)=\exp_{<_{L'}^h}(Q_j)$ for any $j$.
\item for any $L\in \sigma$, the set $\{Q_1,\ldots,Q_r\}$ is the minimal reduced standard basis of $h(I)$ w.r.t. $<_L^h$.
\end{itemize}
\end{cor}
We see that for any cone $\sigma$, there exists a set of $Q_j$ which, for any $L$ in $\sigma$, is the $<_L^h$-minimal reduced Gr\"obner basis of $h(I)$. We shall call this set {\bf the standard basis of $h(I)$ associated with $\sigma$}.

Let us end this paragraph with some remarks and notations.
A cone $\sigma$ of $\EV(h(I))$ is not necessarily open (for example, it may be ``semi open''). Thus here is the definition of $\L(\sigma)$:\\
We first consider the closure $\bar{\sigma}$ of $\sigma$.
There exists $L_1, \ldots,L_q \in \UV$ that we assume to be primitive (with $q\ge 1$ that may be greater than $p$) such that
\[\bar{\sigma}=\{L=r_1 L_1+\cdots r_q L_q; \, r_i \ge 0\}.\]
Assume $q$ to be minimal then $\L(\sigma)$ is the set $\{L_1,\ldots,L_q\}$. The set $\L(\EV(h(I)))$ is then nothing but the union of the $\L(\sigma)$ with $\sigma \in \EV(h(I))$ (remark that it is made of integral elements since the cones $\sigma$ are rational).

We define the interior of $\sigma$ as the set of the (strictly) positive combinations of the $L_i$. We denote by $\rangle L_1,\ldots, L_q \langle$ the open cone generated by the $L_i$. Moreover for $L_1$ and $L_2$ in $\UV$, we denote by $\rangle L_1,L_2 \rangle=\{r_1L_1 +r_2L_2; \, r_1>0, \, r_2\ge 0\}$ the semi open cone that contains $L_2$ and not $L_1$.

In the sequel, we shall write $\EV$ instead of $\EV(h(I))$.

\section{Proof of theorem \ref{theo}}\label{sec:demo1}

Let $\sigma$ be a cone of $\EV$. We have seen in Corollary \ref{cor:V_event_Grob} that there exists a family $Q_1,\ldots,Q_r$ in $h(I)$ which is the minimal reduced standard basis of $h(I)$ w.r.t. $<_L^h$, for any $L\in \sigma$. Now what
happens for a linear form $L$ in the closure of $\sigma$ but not in $\sigma$ (which may hold for an $L$ in $\L(\sigma)$)? Of course the $Q_j$ are not necessarily a standard basis of $h(I)$ w.r.t. $<_L^h$. However, it is possible, and this is the purpose of the following proposition, to construct an order that we will denote by $\lhd_L^\sigma$ (because it will depend on $L$ and $\sigma$) for which the $Q_j$ shall be the minimal reduced standard basis of $h(I)$. Here is the definition of the order in question.\\
First we fix a linear form $L_\sigma$ in the interior of $\sigma$, then for $(\alpha,\mu,\beta,\nu,k)$ and $(\alpha',\mu',\beta',\nu',k')$ in $\N^{n+p+n+p+1}$ we set:
\[(\alpha,\mu,\beta,\nu,k) \lhd_L^\sigma (\alpha',\mu',\beta',\nu',k)
\iff\]
\[\begin{cases}
k+|\beta|+|\nu|<k'+|\beta'|+|\nu'| \\
\text{or } \big( = \text{ and } L(\alpha,\mu,\beta,\nu)<
L(\alpha',\mu',\beta',\nu') \big) \\
\text{or } \big( = \text{ and } = \text{ and } (\alpha,\mu,\beta,\nu)
<_{L_\sigma} (\alpha',\mu',\beta',\nu') \big).
\end{cases}\]

\begin{prop}
Let $L \in \UV$ be in the closure of $\sigma \in \EV$ and let $Q_1,\ldots,Q_r$ be the standard basis of $h(I)$ associated with $\sigma$. Then for any $j=1,\ldots,r$ and for any $L'$ in $\rangle L, L_\sigma \rangle$,
\[\exp_{<_{L'}^h}(Q_j)=\exp_{\lhd_L^\sigma}(Q_j).\]
\end{prop}
As a {\bf consequence} of this proposition, we obtain that {\bf the $Q_j$ form the minimal reduced standard basis of $h(I)$ w.r.t. $\lhd_L^\sigma$}. Indeed, it suffices to apply lemma \ref{lem:utile}

\begin{proof}
For the reader's convenience, figure \ref{fig:Qj_sans_pentes} illustrates the situation of the proposition (with $p=2$, but for higher dimension, we could draw a similar figure by intersecting $\UV$ with the plan generated by $L$ and $L_\sigma$).

\begin{figure}[h!]
\centering
\includegraphics[angle=0, width=11cm]{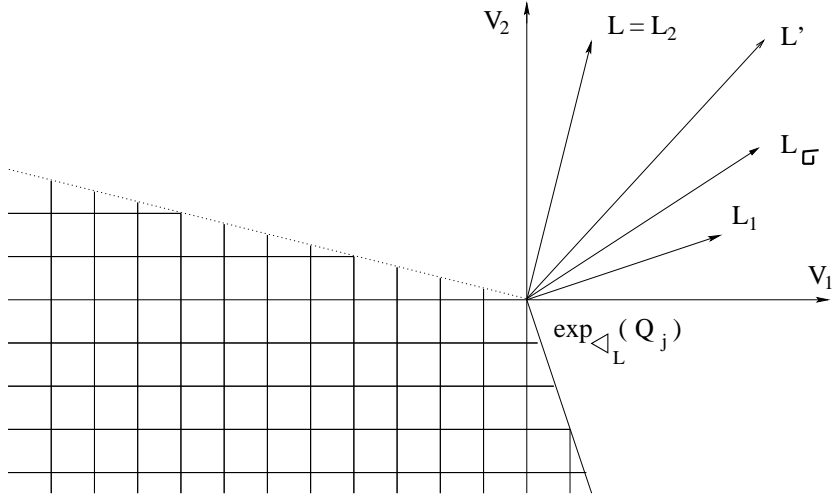}
\caption{Newton diagram of a $Q_j$ associated with $\sigma=\rangle
  L_2, L_1 \rangle$}\label{fig:Qj_sans_pentes}
\end{figure}

\begin{eqnarray*}
\exp_{\lhd_L^\sigma}(Q_j) & = & \exp_{\lhd_L^\sigma}(\sigma^L(Q_j)) \text{ by definition of } \lhd_L^\sigma \text{ and by homogeneity of } Q_j\\
 & = & \exp_{<_{L_\sigma}^h}(\sigma^L(Q_j)) \text{ by definition of
   } \lhd_L^\sigma \text{ et } <_{L_\sigma}^h\\
 & = & \exp_{<_{L_\sigma}^h}(Q_j) \text{ by the following claim}\\
 & = & \exp_{<_{L'}^h}(Q_j) \text{ by corollary \ref{cor:V_event_Grob}}.
\end{eqnarray*}
To end the proof, it remains to proof the following equality.
\begin{aff*}
\[\exp_{<_{L_\sigma}^h}(\sigma^L(Q_j)) = \exp_{<_{L_\sigma}^h}(Q_j).\]
\end{aff*}
Since $\DN(\sigma^L(Q_j)) \subseteq \DN(Q_j)$, we have $\dps \exp_{<_{L_\sigma}^h}(\sigma^L(Q_j)) \le_{L_\sigma}^h \exp_{<_{L_\sigma}^h}(Q_j)$. Then, the reverse inequality suffices to prove the claim.
First let us prove the following:
\[(\star) \qquad \ord^L(Q_j)= L(\exp_{<_{L_\sigma}^h}(Q_j)).\]
For any $L''$ in $\rangle L, L_\sigma \rangle$, we have
\[\ord^{L''}(Q_j)= L''(\exp_{<_{L''}^h}(Q_j))=L''(\exp_{<_{L_\sigma}^h}(Q_j))\]
by definition of $<_{L''}^h$ and by corollary \ref{cor:V_event_Grob}. Thus for any fixed $m$ in $\DN(Q_j)$,
\[L''(\exp_{<_{L_\sigma}^h}(Q_j)) \ge L''(m).\]
Write $L''=(1-\varepsilon) L + \varepsilon L_\sigma$ (with $1 \ge \varepsilon >0$) and make $\varepsilon$ converge to $0$. By continuity, we obtain: $\dps L(\exp_{<_{L_\sigma}^h}(Q_j)) \ge L(m)$, the later holding for any $m$ in $\DN(Q_j)$. Thus the equality $(\star)$ is proved.
As a consequence and by definition of $\sigma^L(Q_j)$, we have
\[\exp_{<_{L_\sigma}^h}(Q_j) \in \DN(\sigma^L(Q_j)).\]
The wanted inequality then follows from the definition of $\exp_{<_{L_\sigma}^h}(\sigma^L(Q_j))$. The claim is proven.
\end{proof}

Recall that $I$ is an ideal in $\Dnp$, $M=\Dnp/I$ and $\delta$ is the class of $1$ in $M$.
Let $\sigma$ be a cone in $\EV$. Denote by $L_1,\ldots,L_q$ the elements of $\mathcal{L}(\sigma)$.
\begin{lem}\label{lem:cle}
Let $i_0 \in \{1,\ldots,q\}$, $m \in V^{L_{i_0}}_{\lambda_{i_0}}(M)$
with $\lambda_{i_0} \in \Q$ and $P \in \Dnp$ such that $P \delta=m$ and
$\ord^{L_{i_0}}(P)>\lambda_{i_0}$ then there exists $P' \in \Dnp$ such that
\begin{itemize}
\item $P-P' \in I$ i.e. $P' \delta=m$,
\item $\ord^{L_{i_0}}(P')<\ord^{L_{i_0}}(P)$,
\item $\ord^{L_i}(P') \le \ord^{L_i}(P)$ for $i\in \{1,\ldots,q\}
  \smallsetminus \{i_0\}$.
\end{itemize}
\end{lem}
In other terms, it is possible to decrease the order w.r.t. one $L_i$ without increasing the order w.r.t. the others $L_i$. Thanks to this lemma we can give a

\begin{proof}[Proof of theorem \ref{theo}]
Let $m \in \bigcap_{L \in \mathcal{L}(\sigma)}
V^L_{L(v)}(M)$ then for $i=1,\ldots,q$, there exists $P_i \in
\Dnp$ such that $P_i \delta=m$ and $\ord^{L_i}(P_i) \le L_i(v)$.
Set $\tilde{P}_1=P_1$. By applying the lemma with $i_0=1$ a finite number of times (the first time with $P=\tilde{P}_1$ and $\lambda_{i_0}=\ord^{L_{i_0}}(\tilde{P}_1)$), we construct $\tilde{P}_2$ such that $\tilde{P}_2\delta=m$ and
$\ord^{L_i}(\tilde{P}_2) \le L_i(v)$ for $i=1,2$. We start again the process with $i_0=3$ and $P=\tilde{P}_2$, etc. After a finite number of steps, we obtain $\tilde{P}_q \in \Dnp$ such that $\tilde{P}_q \delta=m$ and for any $i=1,\ldots,q$, $\ord^{L_i}(\tilde{P}_q) \le L_i(v)$. Finally we can conclude that $m \in
{}^{\sigma}V_v(M)$.
\end{proof}

To end this section, it remains to prove the previous lemma.
\begin{proof}[Proof of lemma \ref{lem:cle}]
To simplify the notations, we do the prove with $i_0=1$ and we set $\lambda=\lambda_{i_0}$ and $\lhd_{L_i}= \lhd_{L_i}^\sigma$. Denote by $Q_1,\ldots,Q_r$ the standard basis of $h(I)$ associated with the cone $\sigma$.

By hypothesis, there exists $P_1 \in \Dnp$ with $P_1 \delta=m$ and $\ord^{L_1}(P_1)\le \lambda$. There exists $l_0,l,l_1 \in \N$ such that $z^{l_0}h(P-P_1)=z^lh(P)-z^{l_1}h(P_1)$. We then set $H=z^lh(P)$, $H_1=z^{l_1}h(P_1)$ and $H_0=H-H_1$ and since $P-P_1$ is in $I$, we have $H_0 \in h(I)$.\\

Consider the division of $H_0$ by the $Q_j$ w.r.t. $\lhd_{L_1}$:
\[H_0=\sum_{j=1}^r q_j Q_j \text{ with }
\mathcal{N}(q_j)+\exp_{\lhd_{L_1}}(Q_j) \subset \Delta_j \text{ for any } j\]
where the $\Delta_j \subset \N^{2n+2p+1}$ form the partition of $\Exp_{\lhd_{L_1}}(h(I))$ associated with the leading exponents of the $Q_j$ (see the division theorem \ref{theo:div}).\\
Since for any $i,j$, the exponents $\exp_{\lhd_{L_1}}(Q_j)$ and $\exp_{\lhd_{L_i}}(Q_j)$ are equal, the previous division is also a division w.r.t. the orders $\lhd_{L_2}, \ldots, \lhd_{L_q}$. As a consequence, for any $i=1,\ldots,q$ and $j=1,\ldots,r$,
\[\ord^{L_i}(H_0) \ge \ord^{L_i}(q_j Q_j).\]
Denote by $J$ the set of the $j\in \{1,\ldots,r\}$ for which $\ord^{L_1}(H_0)=\ord^{L_1}(q_j Q_j)$, then we have
\[\sigma^{L_1}(H_0)=\sigma^{L_1}(H)= \sum_{j\in J} \sigma^{L_1}(q_j)
\sigma^{L_1}(Q_j).\]
Let us consider and denote $\dps W=\sum_{j\in J}
\sigma^{L_1}(q_j)Q_j$. It is an element of $h(I)$.\\
Put $H'=H-W$. We are going to prove the following claims
\begin{enumerate}
\item[(i)] $\ord^{L_1}(H') < \ord^{L_1}(H),$
\item[(ii)] $\ord^{L_i}(H') \le \ord^{L_i}(H)$ for $i=2,\ldots,q$.
\end{enumerate}

\begin{enumerate}
\item[(i)]
Clearly $\sigma^{L_1}(H)=\sigma^{L_1}(W)$. Therefore,
\[H'=(H-\sigma^{L_1}(H)) - (W -\sigma^{L_1}(W)).\]
We can easily see that the two terms in brackets have an $L_1$-order strictly less than that of $H$.

\item[(ii)]
Fix $i$ between $2$ and $q$.\\
By cor. \ref{cor:V_event_Grob}, for any $j$, we have
\begin{equation}\label{ident1}
\exp_{\lhd_{L_i}}(\sigma^{L_1}(Q_j))= \exp_{\lhd_{L_i}}(Q_j).
\end{equation}
On the other hand, by construction of the $q_j$, for any $j$: $\DN(q_j)+ \exp_{\lhd_{L_i}}(Q_j)\subset \Delta_j$, then for any $j$, the following holds
\begin{equation}\label{ident2}
\DN(\sigma^{L_1}(q_j))+\exp_{\lhd_{L_i}}(\sigma^{L_1}(Q_j))\subset
\Delta_j.
\end{equation}
Now, we have $\sigma^{L_1}(W)=\sum_{j\in J} \sigma^{L_1}(q_j)
\sigma^{L_1}(Q_j)$. By relation (\ref{ident2}), we can say that this writing is the result of the division of $\sigma^{L_1}(W)$ by $\{\sigma^{L_1}(Q_j), \, j\in J\}$ w.r.t. $\lhd_{L_i}$, therefore:
\begin{equation}\label{ident3}
\exp_{\lhd_{L_i}}(\sigma^{L_1}(W)) =\max_{j\in J}
\{\exp_{\lhd_{L_i}} (\sigma^{L_1}(q_j) \sigma^{L_1}(Q_j)) \}.
\end{equation}
In a similar way, we can prove
\begin{equation}\label{ident4}
\exp_{\lhd_{L_i}}(W) =\max_{j\in J}
\{\exp_{\lhd_{L_i}} (\sigma^{L_1}(q_j) Q_j) \}.
\end{equation}
As a consequence, thanks to (\ref{ident3}), (\ref{ident4}), (\ref{ident1}), we obtain the equality
$\exp_{\lhd_{L_i}}(\sigma^{L_1}(W)) =\exp_{\lhd_{L_i}}(W)$, which implies in particular $\ord^{L_i}(W)=\ord^{L_i}(\sigma^{L_1}(W))$.
Hence the following relations hold:
\begin{eqnarray*}
\ord^{L_i}(W) & = & \ord^{L_i}(\sigma^{L_1}(W)) \\
 & = & \ord^{L_i}(\sigma^{L_1}(H)) \text{ because }
 \sigma^{L_1}(W)=\sigma^{L_1}(H) \\
 & \le & \ord^{L_i}(H). 
\end{eqnarray*}
Therefore, $\ord^{L_i}(H') \le \ord^{L_i}(H)$.
The two claims are proven.
\end{enumerate}
Now let us specialize $z=1$ (which a morphism between the algebras $\Dnpz$ and $\Dnp$) and set $P'=H'_{|z=1}=P-W_{|z=1}$. Since $W\in h(I)$, we have $W_{|z=1} \in I$ and $P'\delta=m$. After specialization, the claims (i) and (ii) become: $\ord^{L_i}(P') \le \ord^{L_i}(P)$ for any $i=1,\ldots,q$ with a strict inequality for $i=1$. The lemma is proven.
\end{proof}

\section{Proof of theorem \ref{theo2}}\label{sec:demo2}

In this section, we will give the proof of th. \ref{theo2}. We shall firstly state more precisely the theorem in question. We shall see that the proof consists essentially in a refined analysis of the previous lemma \ref{lem:cle}. Indeed, we shall do what we could call a control over the order w.r.t. to the form $V_1$.\\
Recall that in this section $p$ equals $2$.

\begin{nota*}
\
\begin{itemize}
\item
Let $L_1,L_2$ be two non zero forms in $\UV$. Write
$L_i=a_iV_1+b_iV_2$ with $a_i,b_i \ge 0$. We say that $L_1$ is lower
(resp. strictly lower) than $L_2$ if $b_1/a_1
\le b_2/a_2$ (resp. $b_1/a_1 < b_2/a_2$). We denote this notion by
$L_1 \le L_2$ (resp. $L_1 < L_2$). By convention,
$b/0=+\infty$, i.e. any form $L$ is lower than $V_2$.
\item
Let $L_1 \ne L_2$ be in $\UV$ and $H\in \Dnpz$. We say that $H$ is
$L_1$-homogeneous if $H=\sigma^{L_1}(H)$. We say that $H$ is
$(L_1,L_2)$-homogeneous if $H=\sigma^{L_1}(\sigma^{L_2}(H))$.
\item
Let $L$ be a form in $\UV$. We define $\lhd_L$ as the order on
$\N^{n+2+n+2+1}$ given by:
\[(\alpha,\mu,\beta,\nu,k) \lhd_L (\alpha',\mu',\beta',\nu',k) \iff\]
\[\begin{cases}
k+|\beta+\nu|<k'+|\beta'+\nu'| \\
\text{or } \big( = \text{ and } L(\alpha,\mu,\beta,\nu)<
L(\alpha',\mu',\beta',\nu') \big) \\
\text{or } \big( = \text{ and } = \text{ and } (\alpha,\mu,\beta,\nu)
<_{V_1} (\alpha',\mu',\beta',\nu') \big).
\end{cases}\]
Remark that with the notation of the previous sections, if we set $\sigma=\rangle V_1, L \langle$ (with $L\ne V_1$) then $\lhd_L=\lhd_L^\sigma$, and if $L=V_1$ then $\lhd_L=<_{V_1}^h$.
\end{itemize}
\end{nota*}

\begin{figure}[h!]
\centering
\includegraphics[angle=0, width=9cm]{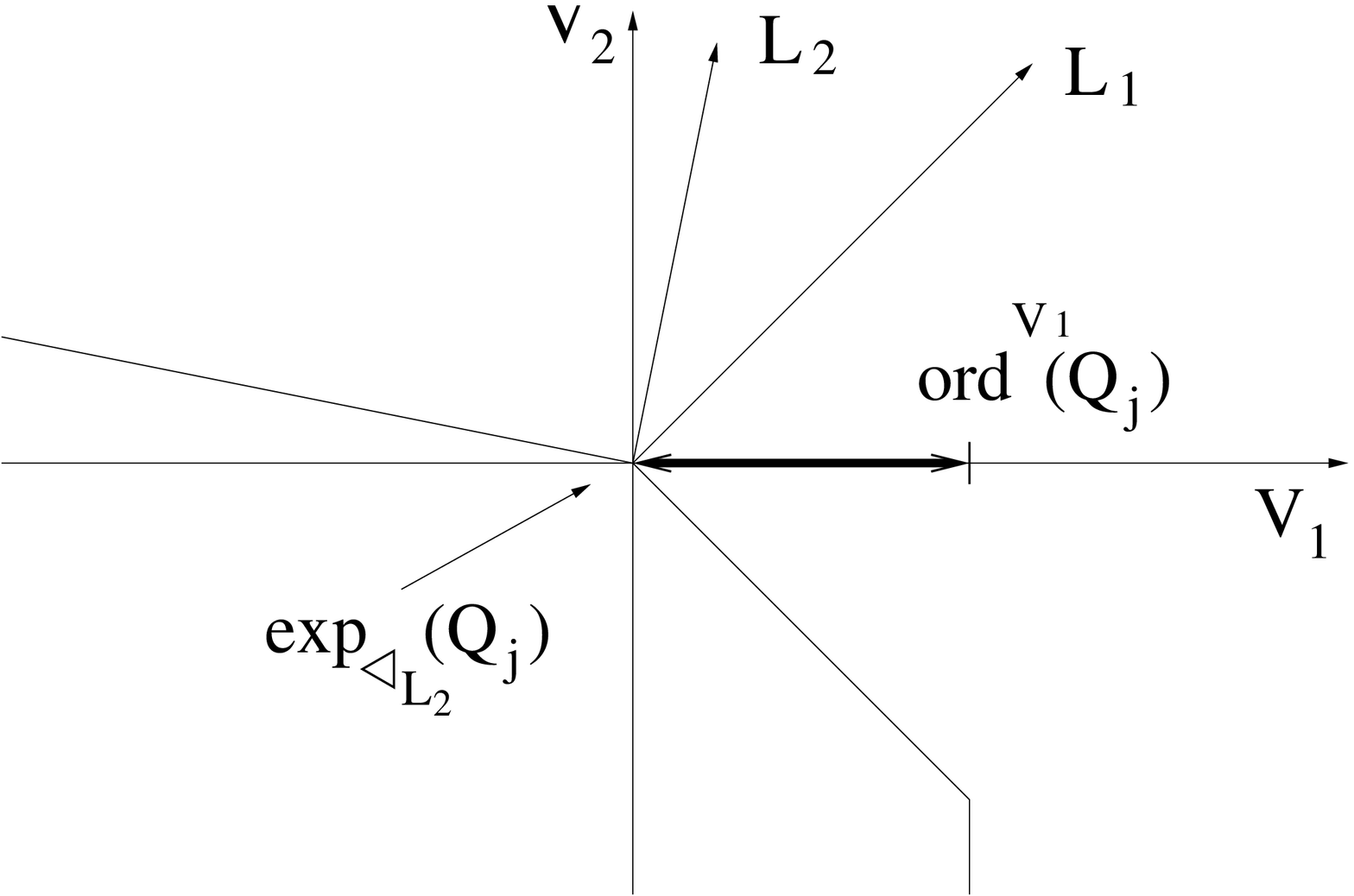}
\caption{$\ord^{V_1}(Q_j)-
  \ord^{V_1}(\exp_{\lhd_{L_2}}(Q_j))$}\label{fig:kappa1}
\end{figure}

Let $\sigma$ be a cone of $\EV$ of (maximal) dimension $2$ and let $\{L_1,L_2\}=\mathcal{L}(\sigma)$ with $L_1 <L_2$. Let $Q_1,\ldots,Q_r$ be the standard basis of $h(I)$ associated with $\sigma$. We define $\kappa_\sigma^1 \in \N$ by:
\[\kappa_\sigma^1=\max\{ \ord^{V_1}(Q_j)-\ord^{V_1}(\sigma^{L_2}(Q_j))
, \, j=1,\dots,r\}.\]
With the previous notations, we have (see figure \ref{fig:kappa1})
\[\ord^{V_1}(\sigma^{L_2}(Q_j))=\ord^{V_1}(\exp_{\lhd_{L_2}}(Q_j)).\]

Now, we define $\kappa^1 \in \N$ as the maximum of all the $\kappa_\sigma^1$ for $\sigma \in \EV$ of maximal dimension. Here a restatement of theorem \ref{theo2}.

\begin{theobis*}
For any $w \in \Z^2$:
\[\Vb_w(M) \subset V_{w+(\kappa^1,0)}(M).\]
\end{theobis*}

\subsection{Control of the $V_1$-order}

Let $\sigma \in \EV$ be a cone of maximal dimension et let $L_1,L_2$ be its primitive generators. Let $m\in \Vb_w(M)$ with $w$ in $\Z^2$, in particular $m\in (V^{L_1}_{L_1(w)}(\Dnd)\delta) \cap (V^{L_2}_{L_2(w)}(\Dnd)\delta) $. Suppose given $P\in \Dnp$ such that $P\delta=m$ and $\ord^{L_1}(P)\le L_1(w)$ and such that $\ord^{L_2}(P) > L_2(w)$. Then we have shown in lemma \ref{lem:cle} how we can construct, in a finite number of steps, an operator $P_\sigma$ such that $\ord^{L_1}(P_\sigma) \le \ord^{L_1}(P)$ (i.e. the $L_1$-order has not increased) and $\ord^{L_2}(P_\sigma) \le L_2(w)$ (i.e. the $L_2$-order has decreased as much as possible). We can wonder what happens concerning the $V_1$-order of $P_\sigma$ compared with that of $P$. We are going to show that this order can be greater but in a controled way. That is the purpose of the next lemma.

\begin{lem}\label{lem:controleV1}
Let $\sigma$ be a cone of maximal dimension in $\EV(h(I))$ and $L_1
\ne L_2$ its primitive generators (which can be out of $\sigma$). Suppose $V_1 \le L_1 < L_2 \le V_2$.\\
Let $w\in \Z^2$ and $m \in V^{L_2}_{L_2(w)}(M)$. Let $P\in \Dnd$ be
such that $P\delta=m$ and $\ord^{L_1}(P)\le L_1(w)$ then we can construct $P_\sigma \in \Dnd$ from $P$ such that:
\begin{description}
\item[(i)] $P_\sigma -P \in I$
\item[(ii)] $P_\sigma \in {}^{\sigma}V_w(\Dnd)$, in particular:
  $\ord^{L_2}(P_\sigma)\le L_2(w)$
\item[(iii)] $\ord^{V_1}(P_\sigma) \le \max\{ \ord^{V_1}(P) \, ,\,
  w_1+\kappa_\sigma^1\}$.
\end{description}
\end{lem}
Statement $(iii)$ justifies the title of this paragraph.
\begin{proof}
If $\ord^{L_2}(P)\le L_2(w)$, it suffices to set $P_\sigma=P$.
Let us then assume that $\ord^{L_2}(P)>L_2(w)$ which implies
$\ord^{V_1}(\sigma^{L_2}(P)) \le w_1$.

By hypothesis, there exists $P_2 \in \Dnd$ satisfying $P_2\delta=m$ and
$\ord^{L_2}(P_2)\le L_2(w)$. We define $H_0=z^{l_0}h(P-P_2)
=z^lh(P)-z^{l_2}h(P_2)$ (there exists integers $l_0$, $l$, and $l_2$
for which such an equality holds), $H=z^lh(P)$ and $H_2=z^{l_2}h(P_2)$.
Let us restart the proof of lemma \ref{lem:cle} with the difference that we work with the form $L_2$ instead of $L_1$. Thus we consider the division of $H_0$ by the standard basis $Q_1,\ldots,Q_r$ w.r.t. the order $\lhd_{L_2}$:\\
$H_0=\sum_{j=1}^r q_j Q_j$ with $\ord^{L_2}(H_0) \ge
\ord^{L_2}(q_jQ_j)$. We denote by $J$ the set of $j$ in $\{1,\ldots,r\}$ for which the later is an equality. Then we set

$W=\sum_{j \in J} \sigma^{L_2}(q_j) Q_j$ et $H'=H-W$.\\
Now we are interesting in comparing $\ord^{V_1}(H)$ and $\ord^{V_1}(H')$.

\begin{affs*}
\
\begin{description}
\item[(a)] $\ord^{V_1}(W) \le w_1+\kappa_\sigma^1$
\item[(b)] $\ord^{V_1}(W)-\ord^{V_1}(\sigma^{L_2}(W)) \le
  \kappa^1_\sigma$
\end{description}
\end{affs*}
Let us prove these claims.\\
{\bf (a):} We have $\ord^{V_1}(\sigma^{L_2}(W))=
\ord^{V_1}(\sigma^{L_2}(H))=\ord^{V_1}(\sigma^{L_2}(P))$ and the later is upper bounded by $w_1$ then if {\bf (b)} is true then so it is for {\bf (a)}.\\
{\bf (b):} As in the previous paragraph, we can show that the
division of $W$ by $\{Q_j,\, j\in J\}$ w.r.t. $\lhd_{L_2}$ is: $W=\sum_{j \in J} \sigma^{L_2}(q_j) Q_j$. Therefore, there exists $j_1\in J$ such that $\exp_{\lhd_{L_2}}(W)=\exp_{\lhd_{L_2}}(m_1 Q_{j_1})$ where $m_1=\mp_{\lhd_{L_2}}(q_{j_1})$ is a monomial of $\sigma^{V_1}(\sigma^{L_2}(q_{j_1}))$. In particular, this implies
\begin{equation}\label{ident5}
\ord^{V_1}(\sigma^{L_2}(W))= \ord^{V_1}(m_1 \sigma^{L_2}(Q_{j_1})).
\end{equation}
Moreover,
\[\ord^{V_1}(W) \le \max\{ \ord^{V_1}(\sigma^{L_2}(q_j) Q_j);\, j \in
J \}.\]
Then let $j_2 \in J$ be such that $\ord^{V_1}(\sigma^{L_2}(q_{j_2})
Q_{j_2})= \max\{ \ord^{V_1}(\sigma^{L_2}(q_j) Q_j),\,  j\in J\}$. By setting $m_2=\mp_{\lhd_{L_2}}(q_{j_2})$, we obtain 
\begin{equation}\label{ident6}
\ord^{V_1}(W)\le \ord^{V_1}(m_2 Q_{j_2}).
\end{equation}
Remark that we may have $j_1=j_2$. However, in any case (see figure \ref{fig:affirmations}), we have:

\begin{figure}[h!]
\centering
\includegraphics[angle=0, width=11cm]{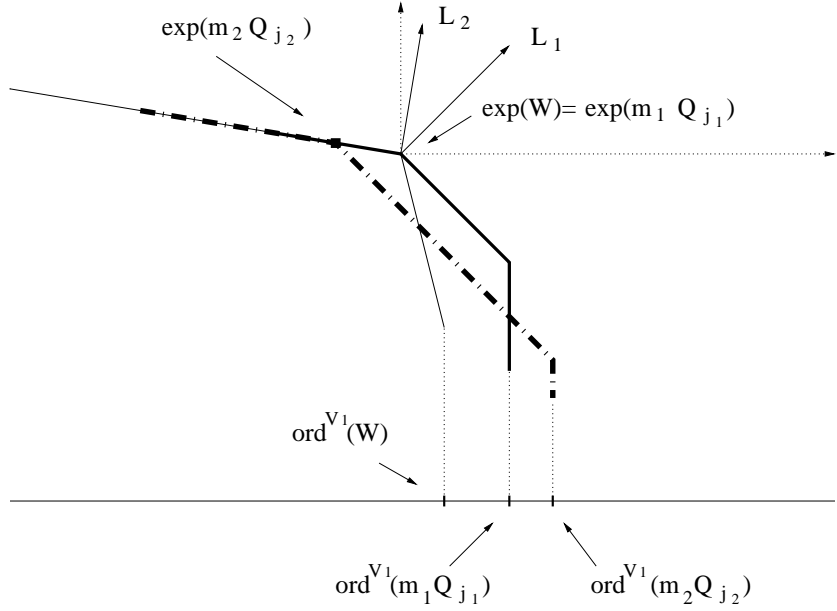}
\caption{Illustration of the claims}\label{fig:affirmations}
\end{figure}

\begin{aff*}
\
\begin{description}
\item[(c)] $\ord^{V_1}(m_2 \sigma^{L_2}(Q_{j_2})) \le \ord^{V_1}(m_1
  \sigma^{L_2}(Q_{j_1}))$.
\end{description}
\end{aff*}
By using this last claim and identities (\ref{ident5}) and (\ref{ident6}), we obtain:
\begin{eqnarray*}
\ord^{V_1}(W) & = & \ord^{V_1}(W) - \ord^{V_1}(m_2 Q_{j_2})\\
 & & + \ord^{V_1}(m_2 Q_{j_2}) -\ord^{V_1}(m_2
 \sigma^{L_2}(Q_{j_2}))\\
 & & +\ord^{V_1}(m_2 \sigma^{L_2}(Q_{j_2}))-
 \ord^{V_1}(m_1 \sigma^{L_2}(Q_{j_1}))\\
 & & +\ord^{V_1}(\sigma^{L_2}(W))\\
 & \le & \ord^{V_1}(m_2 Q_{j_2})
 -\ord^{V_1}(m_2 \sigma^{L_2}(Q_{j_2})) +
 \ord^{V_1}(\sigma^{L_2}(W))\\
 & \le & \kappa^1_\sigma+ \ord^{V_1}(\sigma^{L_2}(W))
\end{eqnarray*}
This proves claim {\bf (b)}. It remains to prove claim {\bf (c)}.\\

The division of $\sigma^{L_2}(W)$ by $\{\sigma^{L_2}(Q_j),\, j\in J\}$
w.r.t. $\lhd_{L_2}$ is written as:
$\sigma^{L_2}(W)= \sum_{j\in J} \sigma^{L_2}(q_j) \sigma^{L_2}(Q_j)$. Therefore,
\[\exp_{\lhd_{L_2}}\big(\sigma^{L_2}(q_{j_2}) \sigma^{L_2}(Q_{j_2})
\big) \unlhd_{L_2} \exp_{\lhd_{L_2}}(\sigma^{L_2}(W)).\]
But
\[\exp_{\lhd_{L_2}}(\sigma^{L_2}(W))= \exp_{\lhd_{L_2}}(W)=
\exp_{\lhd_{L_2}}(m_1\sigma^{L_2}(Q_{j_1}))\]
and
\[\exp_{\lhd_{L_2}}\big(\sigma^{L_2}(q_{j_2}) \sigma^{L_2}(Q_{j_2})
\big)= \exp_{\lhd_{L_2}}(m_2\sigma^{L_2}(Q_{j_2}))\]
then
\[\exp_{\lhd_{L_2}}(m_2\sigma^{L_2}(Q_{j_2})) \unlhd_{L_2}
\exp_{\lhd_{L_2}}(m_1\sigma^{L_2}(Q_{j_1})).\]
But $\ord^{L_2}(m_2\sigma^{L_2}(Q_{j_2}))=
\ord^{L_2}(m_1\sigma^{L_2}(Q_{j_1}))$ then
\[\ord^{V_1}(m_2\sigma^{L_2}(Q_{j_2})) \le
\ord^{V_1}(m_1\sigma^{L_2}(Q_{j_1})).\]
Claim {\bf (c)} is proved.\\
Now let us see how claim {\bf (a)} allows one to prove the third claim of the lemma. We started with $H$ and we have constructed $H'=H-W$.
By {\bf (a)}, we have $\ord^{V_1}(H') \le
\max (\ord^{V_1}(H), w_1+\kappa_\sigma^1)$. We continue the same process with $H'$ instead of $H$, etc. The last element $H_\sigma$ that we obtain satisfies:
$H_\sigma-H \in h(I)$ et $\ord^{V_1}(H_\sigma) \le \max
(\ord^{V_1}(H), w_1+\kappa_\sigma^1)$.\\
We finally set $P_\sigma={H_\sigma}_{|z=1}$, and we have $P_\sigma-P \in
I$ et $\ord^{V_1}(P_\sigma) \le \max (\ord^{V_1}(P),
w_1+\kappa_\sigma^1)$. The lemma is proven.
\end{proof}

\subsection{End of the proof}

\begin{proof}[Proof of theorem 2 bis]
Denote by $L_0=V_1<\cdots <L_q=V_2$ the primitive elements of the $1$-skeleton of $\EV$. For each $i=1,\ldots,q$, denote by
$\sigma_i \in \EV$ the (maximal) cone containing the open cone generated by $L_{i-1}$ and $L_i$.\\
Let $m \in \Vb_w(M)$.\\
Let us prove by an induction on $i$ that for any $i=0,\ldots,q$, there exists $T_i \in \Dnd$ satisfying:
\begin{itemize}
\item $T_i \delta=m$
\item $T_i \in V^{L_i}_{L_i(w)}(M)$
\item $\ord^{V_1}(T_i) \le w_1+\kappa^1$.
\end{itemize}
For $i=0$: $m\in \Vb_w$ then in particular $m \in V^{V_1}_{V_1(w)}$
(note that $V_1(w)=w_1$) then there exists $T_0$ such that $T_0 \delta=m$
and $\ord^{V_1}(T_0)\le w_1 \le w_1+\kappa^1$.\\
Assume the statement to be true at rank $i-1$.\\
Let us apply lemma \ref{lem:controleV1} with $\sigma=\sigma_i$ and $P=T_{i-1}$. Then set $T_i=P_\sigma$ (notations of the lemma). Thanks to the lemma in question, $T_i$ satisfies:
\begin{itemize}
\item $T_i \delta=m$
\item $T_i \in V^{L_i}_{L_i(w)}(M)$
\item $\ord^{V_1}(T_i) \le \max(\ord^{V_1}(T_{i-1}),w_1+\kappa^1)
  = w_1+\kappa^1$
\end{itemize}
Thus, the claim is true for any $i$. The particular case $i=q$ gives us: $m= T_q \delta$, $\ord^{V_2}(T_q)\le w_2$ and $\ord^{V_1}(T_q)
\le w_1+\kappa^1$, which means $m \in V_{w+(\kappa^1,0)}(M)$.
\end{proof}

\begin{rem}
The way we constructed $\kappa$ shows that it is not unique in general. Indeed, by working with $V_2$ instead of $V_1$, we could construct some $\kappa$ of the form $(0,\kappa^2)$.
\end{rem}

\end{document}